\documentclass[preprint,3p,12pt]{elsarticle}

\usepackage{amssymb}

\usepackage{lipsum}
\usepackage{subfigure}
\usepackage{amsfonts}
\usepackage{amsthm}
\usepackage{graphicx}
\usepackage{epstopdf}
\usepackage{amsfonts}
\usepackage{graphicx}
\usepackage{amsopn}
\usepackage{listings}
\usepackage{xcolor}
\usepackage{float}
\usepackage{marvosym}
\usepackage{booktabs}
\usepackage{multirow}
\usepackage{algorithm}
\usepackage{algorithmic}
\usepackage{graphicx}
\usepackage{amsmath}
\usepackage{amsfonts}
\usepackage{amsbsy}
\usepackage{amscd}
\usepackage{amssymb}
\usepackage{latexsym}
\usepackage[mathscr]{eucal}
\usepackage{bm}
\usepackage{threeparttable}
\usepackage{tablefootnote}
\usepackage{threeparttable}

\usepackage{lineno}
\modulolinenumbers[1]

\definecolor{mygray}{gray}{0.9}

\newfont{\bb}{msbm10}

\def\vecx{\mathbf{x}}

\def\vecz{\mathbf{z}}
\def\vecr{\mathbf{r}}
\def\vecf{\mathbf{f}}

\def\vecu{\mathbf{u}}
\def\vecv{\mathbf{v}}

\def\vecs{\mathbf{s}}

\def\matA{\mathbf{A}}
\def\matB{\mathbf{B}}
\def\matC{\mathbf{C}}

\def\matM{\mathbf{M}}

\def\matU{\mathbf{U}}
\def\matV{\mathbf{V}}

\def\matX{\mathbf{X}}
\def\matY{\mathbf{Y}}

\def\matS{\mathbf{S}}

\theoremstyle{proposition}
\newtheorem{proposition}{Proposition}[section]
\newtheorem{definition}{Definition}[section]
\newtheorem{lemma}{Lemma}[section]

\newcommand{\revise}[1]{{\color{blue}#1}}  
\renewcommand{\Re}{\operatorname{Re}}

\journal{arXiv}

\begin{document}

\graphicspath{{figures/}} 


\begin{frontmatter}



\title{A shifted Laplace rational filter for large-scale eigenvalue problems}



\author[a,c]{Biyi Wang}
\ead{wangbiyi20@gscaep.ac.cn}
\author[b]{Karl Meerbergen}
\author[b]{Raf Vandebril}
\author[c,d]{Hengbin An}
\author[c,d]{Zeyao Mo}
            
\affiliation[a]{organization={Graduate School of China Academy of Engineering Physics},
            city={Beijing},
            postcode={100088}, 
            country={China}}

\affiliation[b]{organization={Department of Computer Science, KU Leuven, Celestijnenlaan 200A},
            city={Leuven},
            postcode={3001}, 
            country={Belgium}}

\affiliation[c]{organization={Institute of Applied Physics and Computational Mathematics},
            city={Beijing},
            postcode={100094}, 
            country={China}}

\affiliation[d]{organization={CAEP Software Center for High Performance Numerical Simulation},
            city={Beijing},
            postcode={100088}, 
            country={China}}

\begin{abstract}
We present a rational filter for computing all eigenvalues of a symmetric definite eigenvalue problem lying in an interval on the real axis.
The linear systems arising from the filter embedded in the subspace iteration framework, are solved via a preconditioned Krylov method.

The choice of the poles of the filter is based on two criteria. On the one hand, the filter should enhance the eigenvalues in the interval of interest, which suggests that the poles should be chosen close to or in the interval.
On the other hand, the choice of poles has an important impact on the convergence speed of the iterative method.
For the solution of problems arising from vibrations, the two criteria contradict each other, since fast convergence of the eigensolver requires poles to be in or close to the interval, whereas the iterative linear system solver becomes cheaper when the poles lie further away from the eigenvalues.
In the paper, we propose a selection of poles inspired by the shifted Laplace preconditioner for the Helmholtz equation.

We show numerical experiments from finite element models of vibrations.
We compare the shifted Laplace rational filter  with rational filters based on quadrature rules for contour integration.
\end{abstract}



\begin{keyword}
Generalized eigenvalue problem \sep
Rational filtering \sep
Subspace iteration  \sep
Shifted Laplace

\MSC 65F15 \sep 65N25 \sep 65H17
\end{keyword}

\end{frontmatter}



\section{Introduction}
\label{sect:introd}
Consider the following symmetric definite generalized eigenvalue problem (GEP)
\begin{eqnarray}
\label{eq:GEP}
\matA \vecx = \lambda \matB \vecx, 
\end{eqnarray}
where $\matA \in \mathbb{R}^{n \times n}$ and $\matB \in \mathbb{R}^{n \times n}$
are large, sparse and symmetric matrices with $\matA$ and $\matB$ positive definite.
The GEP~\eqref{eq:GEP} arises from, e.g.,  vibration analysis, quantum mechanics, electronic structure calculations, etc~\cite{saad2011numerical}.
Matrix $\matA$ is the stiffness matrix, $\matB$ is 
the mass matrix, $\lambda$ is an eigenvalue of the matrix pencil $(\matA, \matB)$,
and the nonzero vector $\vecx$ is an associated eigenvector.
The tuple $(\lambda, \vecx)$ is called an eigenpair of the GEP.
The eigenvalues are real and positive.
We are interested in all eigenvalues in a specified interval $(0, \gamma]$ and
their associated eigenvectors.


For large-scale problems, subspace methods
such as shift-and-invert Arnoldi and rational Krylov sequences~\cite{saad2011numerical,ericsson1980spectral,grimes1994shifted,ruhe1984rational,ruhe1994rational_complex_shift}, Jacobi-Davidson~\cite{sleijpen1996jacobi} and residual inverse iteration with projection, also called  generalized Davidson~\cite{morgan1986generalizations}, are commonly used.
Another class of methods is based on a contour integral, utilizing the resolvent operator along a closed contour $\Gamma\subset\mathbb{C}$, which encloses all the desired eigenvalues; this class of methods was proposed in 2003 by Sakurai and Sugiura~\cite{ikegami2010filter,sakurai2003projection}.
Variations of this idea led to a block Krylov subspace method \cite{imakura2014block} and subspace iteration method \cite{beyn2012integral}.
The current paper is closest to the FEAST method~\cite{polizzi2009density}, which is a subspace projection framework. 

In practice, the integration is approximated by using a quadrature rule.
The application of the contour integral requires the solution of a sequence of $\ell$ systems of the form
\begin{equation}
    (\matA - \sigma_j \matB) \vecx = \vecf 
    \quad \text{for} \quad j = 1, \ldots, \ell,
    \label{equ: linear systems for rational function}
\end{equation}
where $\sigma_j$, $j=1,\ldots,\ell$, are the quadrature nodes.
There are several ways to select the quadrature nodes, such as Zolotarev~\cite{guttel2015zolotarev}, midpoint~\cite{xi2016computing}, Gauss-Legendre~\cite{polizzi2009density} and Gauss-Chebyshev~\cite{xi2016computing} rules.
The approximation by a quadrature rule is a rational function of the matrix pencil and is therefore called a \emph{rational filter function}.
See \S\ref{sect: Rational filter} for more details.

The choice of the nodes is not only important for the quality of the filter function, but it has an enormous impact on the performance of the linear system solver.
This paper is about the choice of the nodes so that the filter function satisfies two objectives: (i) good filter quality; (ii) fast linear solution. 
These conditions are not easy to reconcile.
To see this, we have to explain how \eqref{equ: linear systems for rational function} is solved.
A direct method, e.g., based on sparse LU factorization, is usually the first choice, since it is reliable even for ill-conditioned matrices.
In addition, the cost of solving~\eqref{equ: linear systems for rational function} with a direct method is usually independent of $\sigma_j$.
When the problem scale is very large, iterative methods may become preferable.
In this case the cost for solving \eqref{equ: linear systems for rational function} strongly depends on the value of $\sigma_j$.
The main motivation of the paper is therefore to choose $\sigma_j$ to reconcile fast convergence of the eigenvalue solver and convergence of
the iterative linear system solver.


For the Helmholtz equation and the elasticity equation, matrices $\matA$ and $\matB$ are symmetric positive definite, but $\matA - \mu \matB$ with $\mu$ real and positive, is generally indefinite.
Preconditioning an indefinite linear system is considered hard, since the matrix usually is far from an M-matrix.
A particularly interesting idea is the shifted Laplace preconditioner for the Helmholtz equation \cite{erlangga2004class,erlangga2006novel,van2007spectral}, which was inspired on
earlier work
\cite{magolu2000preconditioning,magolu2001incomplete,made2004performance}.
The general idea is not to use a preconditioner for $\matA - \mu \matB$ but for
\begin{equation}
\matC = \matA - \mu (1 \pm \alpha \imath) \matB \quad \text{with} \quad \alpha \neq 0, \mu \in \mathbb{R}^+.
\label{equ: pole away real axis}
\end{equation}
i.e., an imaginary matrix is added. It makes the linear system complex, but easier to solve.

In this paper, we intend to reconcile the two conditions mentioned above by picking the quadrature nodes on the line $\{z = \xi(1 + \alpha \imath), \xi\in\mathbb{R}\}$.
Note that moving the shift away from the spectrum in the shift-and-invert transformation was also suggested by Xi and Saad
\cite{xi2016computing}, Meerbergen and Roose\cite{meerbergen1997restarted}, Ohno, Kuramashi, Sakurai and Tadano\cite{ohno2010quadrature}.
Our analysis shows that a key requirement for this strategy is that the matrix $\matB$ remains well-conditioned.

The outline of this paper is as follows. In Section~\ref{sect:poles away real axis}, we analyze and experimentally show that linear equations with complex shifts are easier to solve. 
In Section~\ref{sect: Rational filter}, we review rational filter methods for eigenvalue computations
and discuss our new strategy of selecting the poles and weights for the rational filter.
We also analyse the influence of the key parameters in the new rational filter.
Furthermore, we compare the filter with the rational filters obtained by using the midpoint, Gauss-Legendre~\cite{polizzi2009density}, and Gauss-Chebyshev~\cite{xi2016computing} rules applied to a circle contour.
In Section~\ref{sect:numer-result}, we present the numerical results of the new rational filter and compare its efficiency with that of the aforementioned filters based on quadrature rules.
From our numerical experiments and theoretical considerations, it follows that the shifted Laplace rational filter does not perform well when the matrix $\matB$ is singular.
We end this paper with the main conclusions in Section~\ref{sect:conclusion}.

The following notation is used throughout this paper.
$\matA \in \mathbb{R}^{n \times n}$ and $\matB \in \mathbb{R}^{n \times n}$
represent $n \times n$ large, sparse and symmetric matrices, $\matA$ and $\matB$ are positive definite.
The $\matB$-inner product of two vectors $\vecu$, $\vecv \in \mathbb{R}^{n}$ is defined by
$
\langle \vecu, \vecv \rangle_{\matB} =
\langle \vecu, \matB \vecv \rangle =
\vecv^{\top} \matB \vecu.
$
The induced norm is $\| \vecv \|_{\matB} = \sqrt{\langle \vecv, \vecv \rangle_{\matB}}$.
The eigenvalues of the matrix pencil $(\matA, \matB)$ are denoted by $\left\{ \lambda_{i} \right\}_{i=1}^{n}$
and are indexed in ascending order, i.e., $0 < \lambda_{1} \leq \lambda_{2} \leq \cdots \leq \lambda_{n}$.
And by $\left\{ \vecx_{i} \right\}_{i=1}^{n}$
we denote the associated eigenvectors as
$\| \vecx_i \|_{\matB} = 1$, $i = 1, 2, \ldots, n$.
For a nonzero vector $\vecz$, the value 
$\theta_{\vecz} = \frac{\vecz^{\top} \matA \vecz}{\vecz^{\top} \matB \vecz}$ is called the
generalized Rayleigh quotient associated with the matrix pencil $(\matA, \matB)$.

\section{Analysis and  experiments for shifted linear systems}
\label{sect:poles away real axis}


In this section, we show the impact of the slope parameter $\alpha$ on the performance of a preconditioned Krylov linear system solver for the linear system
\[
(\matA - \sigma \matB) \vecx = \vecf ,\quad \sigma=\mu(1 + \alpha \imath),
\]
where $\vecf$ is a given right-hand side.
In this section, we give a discussion and numerical evidence to show that it is indeed beneficial to put the poles away from the real axis.
We can state that the larger the slope $|\alpha|$, the further the shift $\sigma$ is away from the real axis.



To illustrate the influence of $\alpha$ on solving linear systems, we use the hollow platform model in \S~\ref{subsect: hollow platform model} as a test example, where the smallest eigenvalue $\lambda_1$ is 104.3308. 
We set the pole $\sigma = \mu (1 + \alpha \imath)$ where $\mu = 104.3$ 
and we use different choices of $\alpha$ from 0 to 1 for the test. 
When $\alpha = 0$, $\sigma$ lies close to the real eigenvalue $\lambda_1$, therefore, the corresponding linear system~\eqref{equ: linear systems for rational function} is ill-conditioned.
We used the iterative solver BiCGStab~\cite{van1992bicgstab} with incomplete LU factorization (with different drop tolerances).
As stopping criteria, we require that the relative residual norm $\|\vecf - (\matA - \sigma \matB) \vecx\|_2/\|\vecf\|_2 < 10^{-10}$ and that the number of iteration is at most 1000. 
The computations were carried out in MATLAB R2024a on a Linux server 
with 56-core Intel processors, and 1 TB of RAM.
We report on the number of iterations and the relative residual norm of the last iteration in Table~\ref{tab: effectiveness of put away poles}.
We point out that BiCGStab  implemented in MATLAB checks the convergence at every half iteration step, thus the number of iterations may not be an integer.

\begin{table}
  \centering
    \begin{tabular}{ccccc}
    Solver & $\mu$ & $\alpha$ & Iteration &  \multicolumn{1}{c}{Residual norm} \\
    \midrule
BiCGStab & 104.3 & 0     & 1000  &  1.60E-01 \\
          &       & 0.05  & 1000  &  4.30E-02 \\
          &       & 0.5   & 1000  &  3.70E-03 \\
          &       & 0.8   & 1000  &  1.50E-03 \\
          &       & 1     & 1000  &  8.20E-04 \\
    \midrule
    BiCGStab+ILU(1e-3) & 104.3 & 0     & 440$^\dagger$   &  3.60E-09 \\
          &       & 0.05  & 342$^\dagger$   &  1.00E-10 \\
          &       & 0.5   & 115   &  9.40E-11 \\
          &       & 0.8   & 86    &  9.60E-11 \\
          &       & 1     & 72    &  5.20E-11 \\
    \midrule
    BiCGStab+ILU(1e-4) & 104.3 & 0     & 60$^\dagger$    &  2.90E-09 \\
          &       & 0.05  & 30.5  &  3.70E-11 \\
          &       & 0.5   & 12    &  2.20E-11 \\
          &       & 0.8   & 9.5   &  3.80E-11 \\
          &       & 1     & 9     &  9.30E-11 \\
    \bottomrule
    \end{tabular}%
\caption{The effect of $\alpha$ on solving the linear systems. The symbol $^\dagger$ indicates that BiCGStab terminated due to stagnation, without reaching the required stop criterion.}
\label{tab: effectiveness of put away poles}%
\end{table}%

From the first block of Table~\ref{tab: effectiveness of put away poles}, we can observe that the linear solver does not converge without preconditioning in 1000 iterations.
As the value of $\alpha$ increases, the residual norm reduces gradually; 
however, it remains very high.
From the second and third block of Table~\ref{tab: effectiveness of put away poles}, we can see that the linear solver converges faster when $\alpha$ gets larger.

\subsection{Analysis about the effect of $\alpha$}
\label{sect: analysis}
\begin{lemma}
    For matrices $\matA, \matB \in \mathbb{C}^{n\times n}$, then 
    \begin{align*}
        \sigma_{max} (\matB\matA) &\leq \sigma_{max} (\matB) \sigma_{max} (\matA), \\
        \sigma_{min} (\matB\matA) &\geq \sigma_{min} (\matB) \sigma_{min} (\matA).        
    \end{align*}
As a consequence, we have for nonsingular $\matB$, that
    \begin{equation}
        \kappa(\matB^{-1} \matA) \leq \frac{1}{\kappa(\matB)}\kappa(\matA).
        \label{equ: condition number inequality}
    \end{equation}  
    \label{lemma}
\end{lemma}

\begin{proposition}
    Consider a symmetric definite matrix pencil $\matA - \lambda \matB$ (with nonsingular $\matB$).
    Given the real scalars $\alpha>0$ and $\mu$.
    We define the pole $\sigma = \mu(1+\alpha \imath)$. The condition number $\kappa(\matC)$ of the matrix $\matC = \matB^{-1}( \matA - \sigma \matB)$ decreases with increasing $\alpha$.
    \label{prop: shifted condition}
\end{proposition}

\proof 
First, define, $\hat{\matC} = \matB^{-1}( \matA - \mu \matB)$.
We define $\lambda_1$ as the eigenvalue of $\matA-\lambda\matB$ closest to $\mu$ and $\lambda_n$ as the eigenvalues furthest from $\mu$.
The condition number $\kappa(\matC)$ is the ratio of the largest and smallest eigenvalues of $\hat\matC$ in absolute value: 
\begin{equation*}
    \kappa(\hat\matC) =  \frac{|\lambda_{n} - \mu|}{|\lambda_{1} - \mu|}.
\end{equation*}

The eigenvalues of matrix $\matC = \matB^{-1}(\matA - \sigma \matB)$ are $\theta_j=(\lambda_j - \mu) - (\alpha \mu) \imath$, $j=1,\ldots,n$.
Since $\lambda_1$ also represents the eigenvalue closest to $\sigma$ in the spectral norm sense, while $\lambda_n$ corresponds the most distant eigenvalue, we have
that the condition number $\kappa(\matC)$ is
\begin{equation*}
    \kappa(\matC) = \frac{\sqrt{(\lambda_n-\mu)^2 + (\alpha \mu)^2}}{\sqrt{(\lambda_1-\mu)^2 + (\alpha \mu)^2}}.
\end{equation*}
Obviously, $\kappa(\matC) < \kappa(\hat{\matC})$ and $\kappa(\matC)$ decreases when $\alpha$ increases.
\qed

\begin{proposition}
    For the matrix pencil $\matA \vecx = \lambda \matB \vecx$ where $\matB$ is not singular, and given a real scalar $\mu$, we define the pole $\sigma = \mu(1+\alpha \imath)$ and we assume $\alpha > 0$ without loss of generality. The condition number $\kappa(\matS)$ of the matrix $\matS = \matA - \sigma \matB$ will get smaller when $\alpha$ gets larger.  
    \label{prop: system condition}
\end{proposition}
\proof
Using the notation from Proposition~\ref{prop: shifted condition} for $\mu, \sigma,  \matC$ and $\hat{\matC}$, and we define $\hat{\matS} = \matA - \mu \matB$.
Lemma~\ref{lemma} then implies that
\[
\kappa(\matS) = \kappa(\matB  \matB^{-1}\matS) \leq \kappa(\matB) \kappa(\matC).
\]
From Proposition~\ref{prop: shifted condition}, we know that
\[
\kappa(\matC) < \kappa(\hat{\matC}),
\]
reapplying Lemma~\ref{lemma} yields
\[
\kappa(\matS) \leq \kappa(\matB)\kappa(\matC) < \kappa(\matB)\kappa(\hat{\matC}) \leq \kappa(\matB) \kappa(\hat{\matS}) \frac{1}{\kappa(\matB)} = \kappa(
\hat{\matS}).
\]
Therefore, we know $\kappa(\matS) < \kappa(\hat{\matS})$ and the condition number $\kappa(\matS)$ will get smaller when $\alpha$ gets larger.
\qed

We can conclude that putting the pole $\sigma = \mu (1 + \alpha \imath )$ away from the real axis makes the linear system easier to solve, as confirmed by the results in Table~\ref{tab: effectiveness of put away poles}.

\subsection{Case of singular $\matB$}
\label{sect: singular case}
In this subsection, we numerically demonstrate that placing poles away provides no benefits when the matrix $\matB$ is singular.
The reason is that increasing $\alpha$ makes $\matA-\sigma\matB$ more similar to the singular matrix $\matB$, and therefore the linear system is hard to solve.

We use the matrix pencil BCSSTK38 and BCSSTM38~\cite{matrixcollection}, which is obtained from the FEM discretization of Boeing airplane engine component model, as a test example.

For this case, the smallest eigenvalue $\lambda_1$ is $2.0419\times 10^4$. 
We set the pole $\sigma = \mu (1 + \alpha \imath)$ where $\mu = 2\times 10^4$, the relative gap between $\mu$ and $\lambda_1$ is $(\lambda_1 - \mu) / (\lambda_1 + \mu)$ = 0.0104.
We use different choices of $\alpha$ from 0 to 10 for the test. 

\begin{table}[htbp]
  \centering
  \caption{The results for solving the associated linear systems with the singular matrix $\matB$.}
    \begin{tabular}{ccccc}
    Solver & $\mu$ & $\alpha$ & Iteration  & \multicolumn{1}{c}{Residual norm} \\
    \midrule
    BiCGStab & 2.00E+04 & 0     & 1000  &  6.80E-01 \\
          &       & 0.5   & 1000  & 7.00E-01 \\
          &       & 1     & 1000  &  6.60E-01 \\
          &       & 10    & 1000  &  7.90E-01 \\
    \midrule
    BiCGStab+ILU(1e-4) & 2.00E+04 & 0     & 1000  &  1.00E+00 \\
          &       & 0.5   & 1000  &  1.00E+00 \\
          &       & 1     & 1000  &  1.00E+00 \\
          &       & 10    & 1000  &  1.00E+00 \\
    \bottomrule
    \end{tabular}%
  \label{tab: singularB}%
\end{table}%

From Table~\ref{tab: singularB}, we can observe that the linear solver can not converge in 1000 iterations, even when $\sigma$ is far away from the real axis.
We can also see that the relative residual norm does not differ much for the different $\alpha$.
\section{Rational filter method}
\label{sect: Rational filter}

Rational filter methods to solve eigenvalue problems often arise from a contour integral method~\cite{sakurai2003projection,polizzi2009density,beyn2012integral}.
A contour integral method uses the resolvent operator along a closed contour $\Gamma\subset\mathbb{C}$, which encloses all the desired eigenvalues.
Here,
\begin{equation}
    \matM = \frac{1}{2\pi i} \int_{\Gamma} (\matA - \sigma \matB)^{-1} \matB d\sigma, \quad p \geq 0,
    \label{equ: moment matrix}
\end{equation}
is a matrix whose action on a vector filters away the components of invariant subspaces associated with the eigenvalues outside the contour.
In practice, exact integration is approximated by applying a quadrature rule, and the approximation of~\eqref{equ: moment matrix} will be 
\begin{equation*}
   \matM \approx \Phi(\matA, \matB) = \frac{1}{2\pi i} \sum^{\ell}_{j=1} w_j (\matA - \sigma_j \matB)^{-1} \matB, 
\end{equation*}
where $\sigma_1,\ldots,\sigma_\ell$ are the quadrature nodes along the domain boundary $\Gamma$ and $w_j$ are the corresponding quadrature weights.
The name rational filter arises from the filter properties when $\Phi(\matA, \matB)$ is applied on a vector $\vecv$.
Let $(\lambda_i,\vecx_i)$, $i=1,\ldots,n$, be the eigenpairs of $\matA-\lambda \matB$.
Decompose $\vecv = \sum_{i=1}^{n} \xi_i \vecx_i$, then 
\begin{equation*}
\Phi(\matA, \matB) \vecv = \sum_{i=1}^{n} \xi_i \Phi (\lambda_i) \vecx_i,
\end{equation*}
where $\lambda_i$ is mapped (or filtered) by the rational function
\begin{equation}
\Phi (\lambda) = \sum_{j=1}^{\ell} \frac{w_j}{\lambda - \sigma_j}.
\label{equ: rational function filter}
\end{equation}
We discuss further the choice of nodes and weights.

Since $\matA$ and $\matB$ are real symmetric and the domain of interest is an interval on the real axis, we choose the $\ell$ nodes as complex conjugate points in the complex plane. Let $\ell=2N$, then, with a proper ordering of the nodes, we have $\sigma_{N+j}=\overline{\sigma_j}$ and $w_{N+j}=\overline{w_j}$, $j=1,\ldots,N$.
As a result,
\[
w_j (\matA - \sigma_j \matB)^{-1} \matB + \overline{w_j } (\matA - \overline{\sigma_j} \matB)^{-1} \matB
 = 2 \mbox{Re}(w_j (\matA - \sigma_j \matB)^{-1} \matB),
\]
which implies that $N$ complex valued linear systems have to be solved when the filter is applied to a real vector.

An overview of the rational filter method for computing eigenpairs of the symmetric definite matrix pencil~\eqref{eq:GEP} is given in Algorithm~\ref{alg:Rational filter method}.
The method aims to find $L$ eigenvalues and associated eigenvectors in the interval $(0,\gamma]$.
The matrix $\matU$ represents the filtered vectors.
The projected matrix pencil $\hat{\matA} = \matU^{\top} \matA \matU$ and $\hat{\matB} = \matU^{\top} \matB \matU$ is formed.
The approximated eigenpairs are extracted from the subspace projection on $\matA$ and $\matB$, by using Rayleigh-Ritz method.
The stopping criterion is $\|\matA \vecx_i - \theta_i \matB \vecx_i\|_2/|\theta_i| \|\matB\vecx_i\|_2 \leq \tau$ for $i=1,\ldots,L$, with $\tau$ a prescribed tolerance.
\begin{algorithm}[htbp]
	\caption{Rational filter method for GEP~\eqref{eq:GEP}}
	\label{alg:Rational filter method}
	\begin{algorithmic}[1]
		\REQUIRE $\matA \in \mathbb{R}^{n\times n}$, $\matB \in \mathbb{R}^{n\times n}$,  interval $(0, \gamma]$, 
        poles $\sigma_1,\ldots,\sigma_N \in \mathbb{C}$ and corresponding weights $w_1,\ldots,w_N\in \mathbb{C}$, random vectors $\matV_1 \in \mathbb{R}^{n\times L}$.
        \\
		\ENSURE converged eigenpairs $[\Theta, \matX]$.\\
		\FOR{$k = 1,2,\ldots$ until all of the wanted eigenpairs have converged}
                \FOR{$j = 1,2,\ldots$, $N$}
                    \STATE Solve $\matY_j = (\matA - \sigma_j \matB)^{-1}\matB\matV_k$.
                \ENDFOR
    \STATE Compute $\matU = \sum_{j=1}^N 2\Re(w_j \matY_j)$.
                \STATE Form $\hat{\matA} = \matU^\top \matA \matU$ and $\hat{\matB} = \matU^\top \matB \matU$.
                \STATE Solve the eigenvalue problem $\hat{\matA} \vecs_i = \theta_i \hat{\matB} \vecs_i$ and let $\vecx_i=\matU \vecs_i$, with $\theta_i$ sorted in ascending order.
                \STATE Compute the relative residual norms $\rho_i = \frac{\|\matA \vecx_i - \theta_i \matB \vecx_i\|_2}{|\theta_i| \|\matB\vecx_i\|_2}$ for $i=1,\ldots,L$.
                \STATE Let $\matV_{k+1}=[\vecx_1,\ldots,\vecx_L]$. 
            \ENDFOR
	\end{algorithmic}
\end{algorithm}



The approximation cost depends on the number of quadrature nodes (or poles), which should therefore be kept as low as possible.
That is, the quality of the filter should be compared with the cost of applying the filter.
The main difference between the various rational filter methods is the selection of poles $\sigma_j$, and corresponding weights $w_j$.
Austin~\cite{austin2015computing} selects poles as Chebyshev roots and the corresponding barycentric weights are used.
The main advantage of this method is that poles and arithmetic are real valued.
However, a potential danger is that a pole accidentally lies close to an eigenvalue.
This blows up the filter and enhances this particular eigenvalue and filters all other eigenvalues.
Xi~\cite{xi2016computing} uses fixed poles and the weights are computed from a least-square problem.
Xi employs a filter with poles of higher multiplicity to enable Krylov subspace recycling, improving the efficiency of filter applications.
In the latest FEAST method, the selection of poles and weights is based on Zolotarev's quadrature rule~\cite{guttel2015zolotarev}. Also here, the efficiency of this strategy depends on the position of the poles.


\subsection{The choice of poles}
\label{subsec: optimal poles}

As an alternative to quadrature rules, Ohno~\cite{ohno2010quadrature} proposed to use two horizontal lines away from the real axis as the integral path for approximating the spectral projection operator. 
Xi~\cite{xi2016computing} selected some fixed poles, which are farther away from the real axis than the quadrature poles.
The weights are obtained from an optimization problem that aims for the best approximation of the step function $h: h(z)=1 \text{ for } z\in[-1,1], \text{ and } 0 \text{ elsewhere}$.

In this paper, we propose a different strategy, that we call the Shifted Laplace Rational Filter (SLRF): it reconciles filter quality and the use of the shifted Laplace preconditioner.
The poles and weights are obtained from an optimization problem that aims for the best approximation of the step function $h$ on the positive real axis.
In this way, we aim to amplify the function values in the interval $(0,\gamma]$ and damp the values in $(\gamma, \infty)$.
We do not need to take into account approximations of the step function outside the positive real axis, since there are only positive real eigenvalues.
The poles lie on the lines
\begin{equation}
    \mathcal{L}_{\pm}: y = x(1 \pm \alpha \imath), \quad x \in (0,\gamma] \ \text{and} \ \alpha > 0,
    \label{equ: fixed line}
\end{equation}
where $\alpha$ is a user-defined slope.
We pick poles $\sigma_j = x_j(1 + \alpha \imath)$, $j=1,\ldots,N$.
As the poles and weights are complex conjugate, we construct the filter
\begin{align}
    \Phi(x) & 
    = \sum_{j=1}^{N}\frac{w_j}{x - x_j(1 + \alpha \imath)} + \frac{\overline{w_j}}{x - x_j(1 - \alpha \imath)}, \quad  x_j \in (0,\gamma]\ \text{and} \ \alpha > 0.
    \label{eqn: rational filter}
\end{align}
The parameters $x_j, w_j$ for $j=1,\ldots,N$ are computed from the optimization problem
\begin{equation*}
    \min_{x_j, w_j} \Vert (\Phi - h)\Vert_{\delta}^2 ,
    \label{equ: optimize problem}
\end{equation*}
where $h$ is the step function which has value one for $x \in [0,\gamma]$ and zero for $x$ elsewhere.
The $\delta$ inner product is inspired by \cite{xi2016computing}:
\begin{equation*}
    <f, g> = \int_{-\infty}^{\infty} \delta(t) f(t) \overline{g(t)} dt,
\end{equation*}
where $\delta(t)$ has the following distribution
\begin{equation*}
    \delta(t) = \left\{ 
                    \begin{array}{ll}
                        0 & \text{if} \quad t > \kappa \ \text{or} \ t < 0, \\
                        \beta & \text{if} \quad 0 \leq t \leq \gamma, \\
                        1 & \text{elsewhere}.
                    \end{array}
               \right.
\end{equation*} 
We note that $\beta$ is a relaxation parameter. With a smaller $\beta$, the difference between $\Phi$ and $h$ in $(0,\gamma]$ can be larger, resulting in a larger slope at $x=\gamma$.
We choose $\kappa=10\gamma$, and $\beta$ is set to either 1 or 0.01 as suggested in~\cite{xi2016computing}.
We used the trapezoidal rule to compute the inner product, and solved the optimization problem.

To further compare different filters, we use the \textit{separation factor} of the filter function, as defined in~\cite{xi2016computing}. This factor is used to characterize the quality of the filter function, in the context of subspace iteration.
\begin{definition}
    Let $\Phi(x)$ be a filter function, we call the separation factor of $\Phi$ its absolute value of the derivative at $x = \gamma:$ $|\Phi'(\gamma)|$. See~\rm{\cite{xi2016computing}}.
    \label{def: separation factor}
\end{definition}

For eigenvalue problems, a larger value of $|\Phi'(\gamma)|$ indicates better separation between the wanted eigenvalues and the unwanted ones. 
With subspace projection, the convergence rate for the desired eigenpairs is faster with a larger $|\Phi'(\gamma)|$~\cite{saad2011numerical,guttel2015zolotarev,xi2016computing,gavin2018krylov}.

\begin{figure}[htbp]
\hspace{-1em}
\includegraphics[scale=0.5]{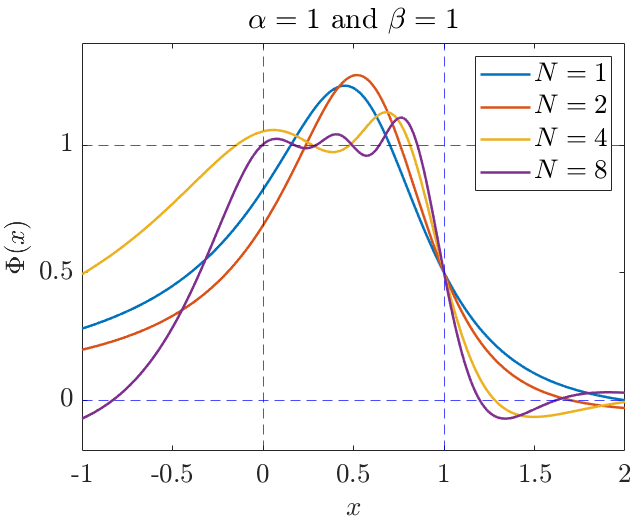}
\hspace{0.5em}
\includegraphics[scale=0.5]{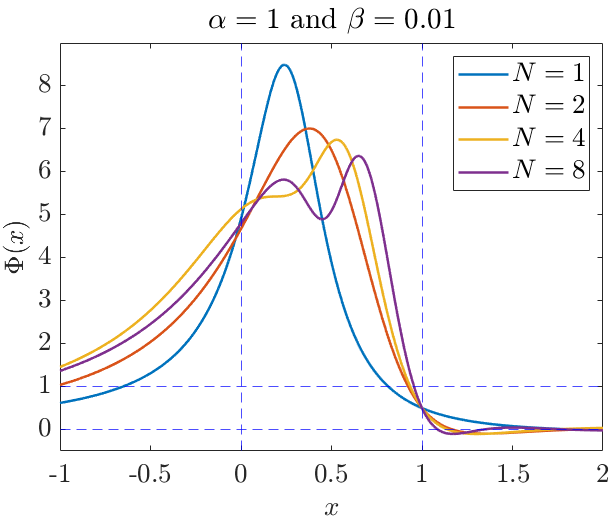}
\caption{The rational filters $\Phi(x)$ are plotted  over the interval $[-1, 2]$ for $\gamma=1$ and $\alpha = 1$, where $\beta = 1$ (left) and $\beta = 0.01$ (right) for different number of poles $N = 1, 2, 4, 8$.}
\label{fig: rational filters}
\end{figure}

The SLRF has three user-defined parameters: the slope $\alpha$ of the straight line, the number of poles $N$ and the value of $\beta$ associated to the inner product.
We now discuss the choice of the parameters.
Figure~\ref{fig: rational filters} shows the influence of the number of poles $N$ on the quality of the rational filter, for fixed $\alpha$.
It is no surprise that the \textit{separation factor} $|\Phi'(\gamma)|$ becomes larger as the number of poles $N$ increases.

We observe that when $\beta = 0.01$, the filter performs better: a smaller value of $\beta$ relaxes the requirement
to approximate the filter value one in the interval $(0, \gamma]$. 
As a result, the filter has larger values in $(0, \gamma]$, smaller values outside $(0, \gamma]$, and a
sharper slope at the endpoint $x = \gamma$, compared to the case when $\beta = 1$.

\begin{figure}[htbp]
\centering
\includegraphics[scale=0.5]{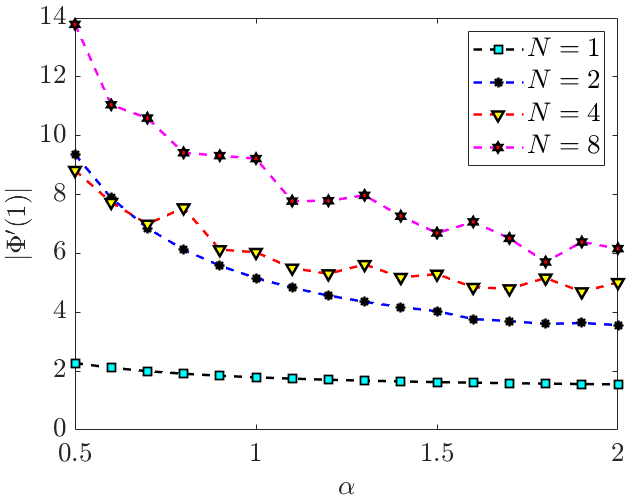}
\caption{Comparison of \textit{separation factors} $|\Phi'(1)|$ of the rational filter using $N = 1, 2, 4, 8$ poles and different slope $\alpha$.}
\label{fig: effect_of_alpha}
\end{figure}

To illustrate the influence of slope $\alpha$ on the \textit{separation factor}, we plot the \textit{separation factor} as a function of $\alpha$ for $\alpha\in[0.5, 2]$ in Figure~\ref{fig: effect_of_alpha}.
We do not show results for $\alpha < 0.5$, because the filter function has high peaks near the real parts of the poles, $x = x_j, j = 1, 2, \ldots, N$.
From Figure~\ref{fig: effect_of_alpha}, we deduce that increasing $\alpha$, not surprisingly, leads to a reduction in \textit{separation factor}.

\begin{figure}[htbp]
\hspace{-1em}
\includegraphics[scale=0.42]{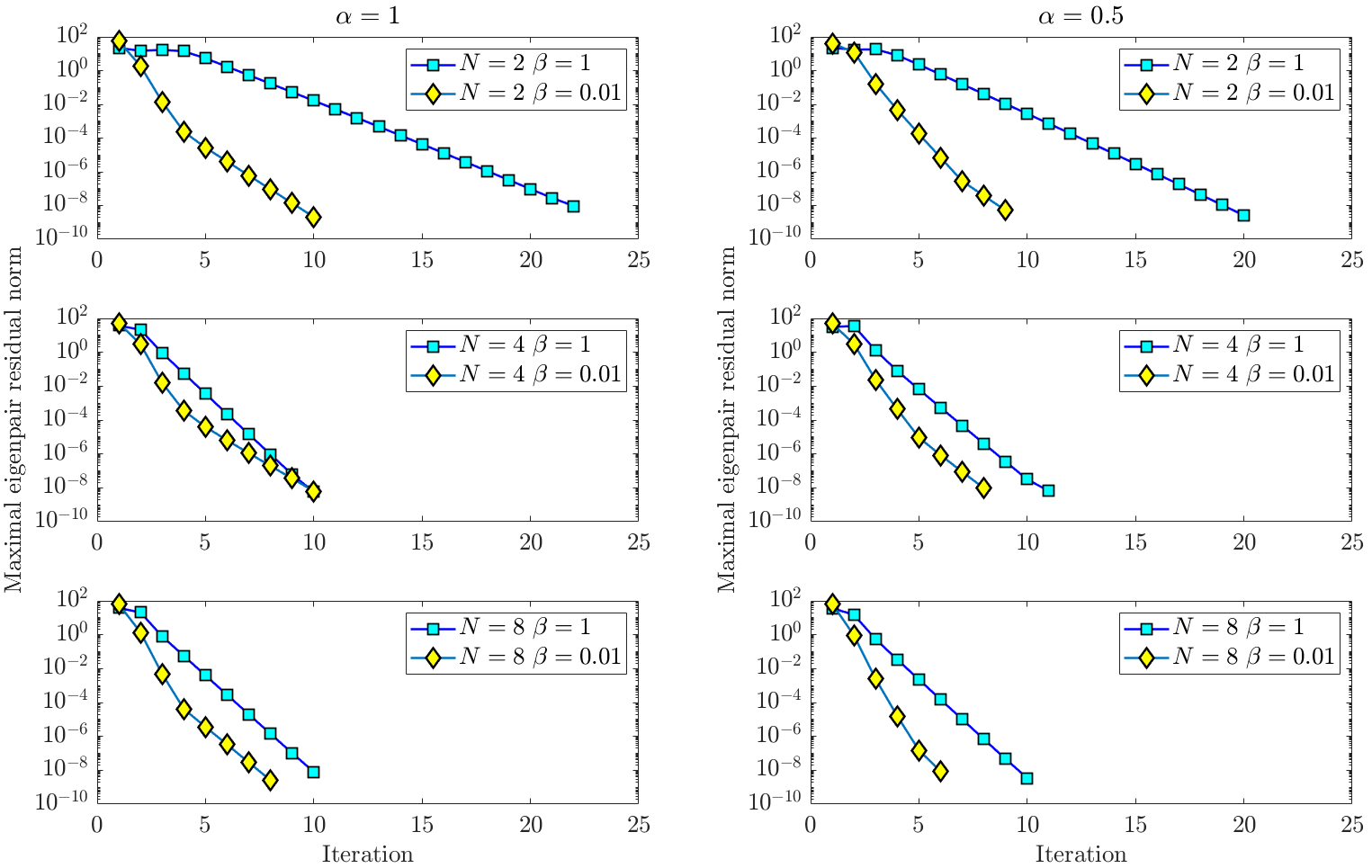}
\caption{The convergence curves of the maximum relative residual norm for the 2D beam model are presented for different filter parameters. The results are shown for $N = 2$ (top), $N = 4$ (middle), $N = 8$ (bottom), with variations in the slope: $\alpha = 1$ (left) and $\alpha = 0.5$ (right), and filter parameters $\gamma=1$, $\beta = 1$ and $\beta = 0.01$ in each subfigure. The maximum residual norm is reported for approximate eigenvalues located within the interval (0, 337.4505], which contains the first 100 smallest eigenvalues, with a search space dimension of $L=120$. The associated linear systems are solved by a complete LU factorization}.
\label{fig: parameters_comp}
\end{figure}

In Figure~\ref{fig: parameters_comp}, we present a detailed numerical comparison of the convergence rate of eigenpairs by using different rational filter parameters, when solving the eigenvalues of the 2D beam model (will be represented in \S\ref{subsect: beam model}).
In each sub-figure, the vertical axis represents the maximum relative residual norm of eigenpairs corresponding to the eigenvalues within the specified interval, and the horizontal axis represents the number of outer iteration. 
We show the results for $\beta = 1$ and $\beta = 0.01$.


In Figure~\ref{fig: parameters_comp}, the first column shows the results when $\alpha = 1$, and the second column shows the results for the smaller $\alpha = 0.5$.
From top to bottom, the number of poles $N$ increases.
From the sub-figure in the upper left corner to the sub-figure in the lower right corner,
we can see that the convergence rate is faster as the number of poles $N$ increases and the slope $\alpha$ gets smaller, particularly for the case $\beta = 1$.

\subsection{Comparison with other filters based on quadrature rules}
\label{subsec: comparison}

We compare the Shifted Laplace Rational Filter with other filters based on classical quadrature rules: midpoint~\cite{xi2016computing}, Gauss-Legendre~\cite{polizzi2009density} and Gauss-Chebyshev quadrature rules~\cite{xi2016computing}, each with $N=3$ poles (in the upper plane).
For SLRF, we set slope  $\alpha = 0.5$ and $\beta = 0.01$.

\begin{figure}[htbp]
\hspace{-0.5em}
\includegraphics[scale=0.42]{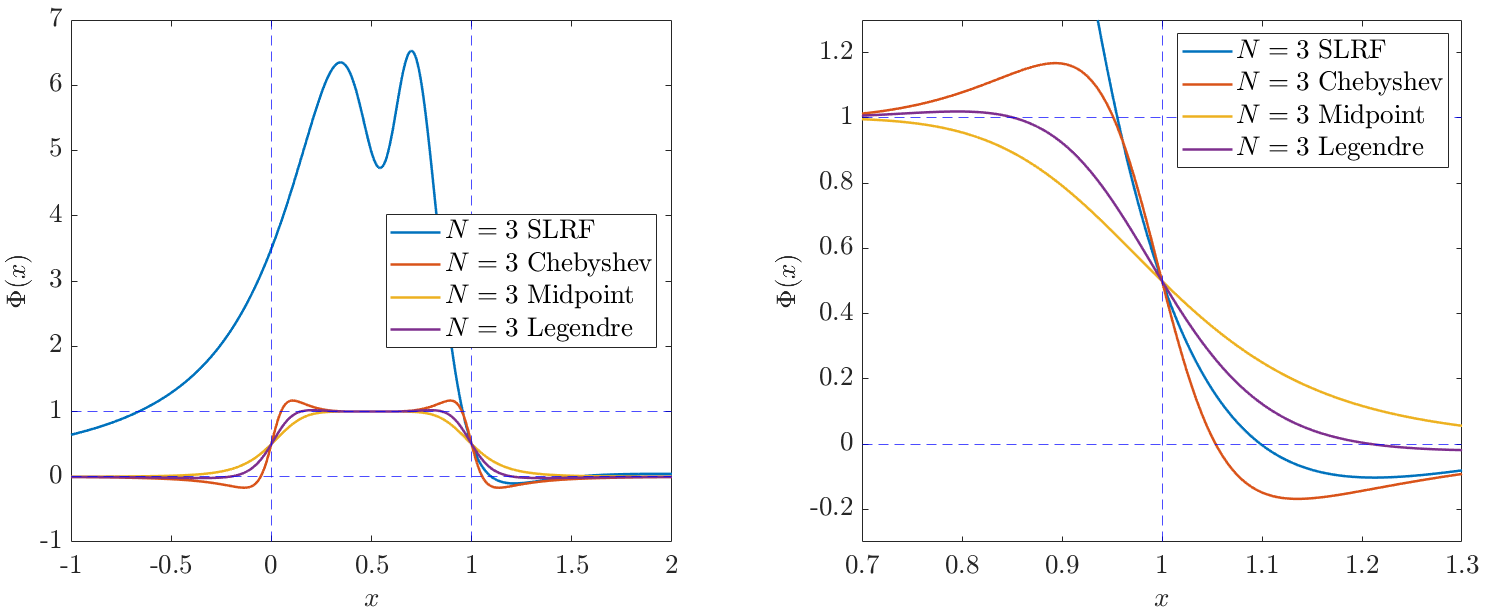}
\caption{Left: the rational filter $\Phi(x)$ for $\gamma=1$ over the interval $x\in[-1, 2]$ for $N = 3$; right: a zoomed-in view at $x=1$.}
\label{fig: comparison}
\end{figure}

From Figure~\ref{fig: comparison} we can see that the \textit{separation factor} of SLRF is larger than that of the midpoint filter and the Gauss-Legendre filter.
The Gauss-Chebyshev filter has the sharpest slope at the endpoint $x = \gamma=1$.
At first sight, there does not seem to be an advantage of SLRF, but this appears
when the linear systems are solved inexactly, as we now show.

The computational task is solving the eigenvalues located in interval $(0, 58.2570]$, which encloses the first 20 smallest eigenvalues of the 2D beam model in \S~\ref{subsect: beam model}.
In Table~\ref{tab: average_iter},  we report the average number of iterations to solve the linear equations associated with each pole for the aforementioned filters.
The detailed setting of algorithm and linear solver can be found in \S~\ref{subsect:numerical-results}.



\begin{figure}[htbp]
\centering
\includegraphics[scale=0.42]{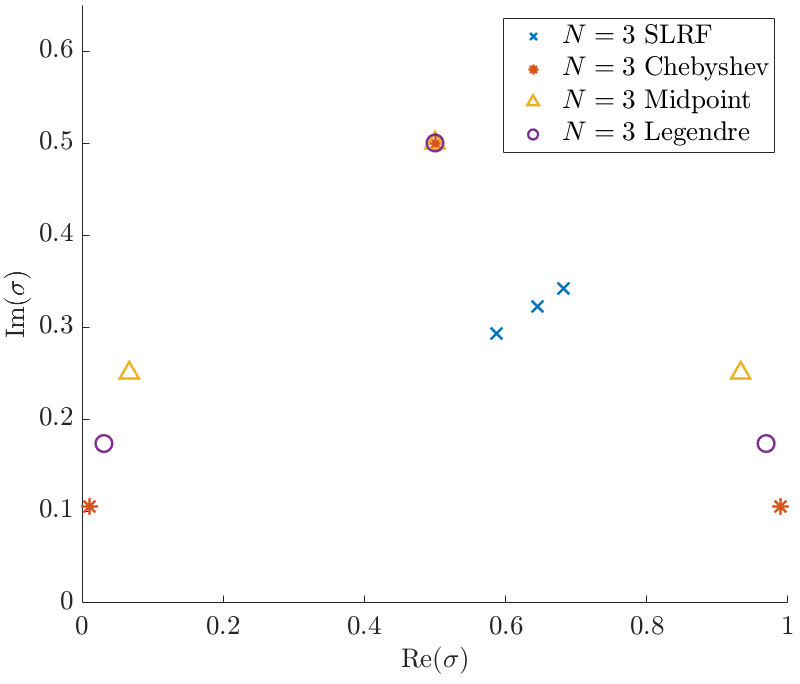}
\caption{The corresponding poles position of four rational filters in the upper half plane.}
\label{fig: N3poles}
\end{figure}

\begin{table}[htbp]
  \centering
  \caption{The average iteration number required to solve the associated linear equations for each pole defined by different filters.}
    \begin{tabular}{cccccc}
    \toprule
    \multirow{2}[2]{*}{Model} & \multirow{2}[2]{*}{Interval} & \multirow{2}[2]{*}{Filter} & \multicolumn{3}{c}{Iter$_{\rm avg}$} \\
          &       &       & $\sigma_1$    & $\sigma_2$    & $\sigma_3$ \\
    \midrule
    \multirow{4}[2]{*}{Beam} & \multirow{4}[2]{*}{(0, 58.2570]} & Midpoint & 25.8  & 27.4  & 58.3 \\
          &       & Gauss-Legendre & 26.3  & 27.4  & 78.1 \\
          &       & Gauss-Chebyshev & 26.5  & 27.4  & 96.1 \\
          &       & SLRF  & 40.7  & 39.5  & 39.2 \\
    \bottomrule
    \end{tabular}%
  \label{tab: average_iter}%
\end{table}%


We also plot the position of the poles of these rational filters, in Figure~\ref{fig: N3poles}. The poles $\sigma_1,\sigma_2,\sigma_3$ are ordered by increasing real part.
From Figure~\ref{fig: N3poles} and Table~\ref{tab: average_iter}, we observe that if the pole is closer to the real axis (i.e., smaller $\alpha$), the corresponding average number of linear iteration Iter$_{\rm avg}$ is larger.
A typical example is the pole $\sigma_3$ defined by the Gauss-Chebyshev filter, which is the closest to the real axis, and incurs the highest computational cost for solving the associated linear systems.
The pole $\sigma_2$ is identical for all three quadrature rules as shown in Figure~\ref{fig: N3poles}, so the cost of solving the linear systems associated with $\sigma_2$ remains the same.
As a consequence, the location of poles can indicate the difficulty of solving the associated linear systems~\eqref{equ: linear systems for rational function}.





\section{Numerical experiments}
\label{sect:numer-result}
\subsection{Finite element model of vibration}
\label{subsect: eigen-problem}
We consider the following vibration system
\begin{eqnarray}
\label{eq:governing equation}
\left\{
\begin{array}{rcl}
\bm{\sigma} \cdot \nabla -\rho \ddot{\mathbf{u}} &=& 0, \quad \text{in} \; \Omega, \\
\bm{\sigma} \cdot \mathbf{n} = \mathbf{\overline{t}}  &=& 0, \quad \text{on} \; \Gamma_N, \\
\mathbf{u}                                       &=& 0, \quad \text{on} \; \Gamma_D,
\end{array}
\right.
\end{eqnarray}
where $\mathbf{u}$ is the displacement, 
$\bm{\sigma}$ is the stress, and 
$\rho$ is the density of the material.
The domain is denoted by $\Omega$, 
$\Gamma_D$ denotes the Dirichlet boundary condition and 
$\Gamma_N$ denotes the Neumann boundary condition, 
$\mathbf{n}$ is the outer normal vector of the boundary of the domain $\partial \Omega$.
The derivation of the weak formulation of~\eqref{eq:governing equation}, which is useful for  finite element discretization, can be found in the Appendix of~\cite{wang2025new}.

In the following, we introduce three mechanical models for numerical tests.
These models have the same material parameters, 
the Young's modulus is $21.5\text{N}/\text{m}^2$,
the Poisson's ratio is 0.29,
the density of the material is $1.0\text{kg}/\text{m}^2$.
The corresponding discretized matrix pencils are generated by using the open source software FreeFEM++~\cite{FreeFem}.



\subsubsection{2D beam model}
\label{subsect: beam model}
\begin{figure}[htbp]
\centering
\includegraphics[scale=0.22]{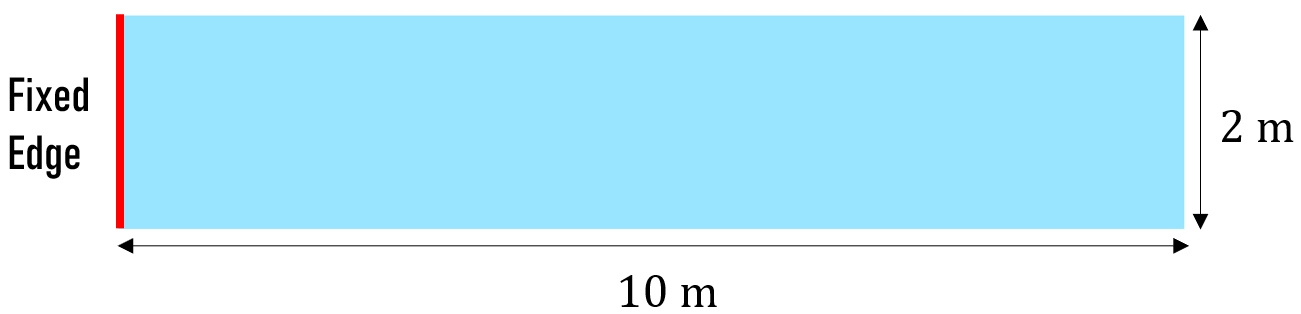}
\caption{The shape of 2D beam model.}
\label{model:beam}
\end{figure}

The first model is a two dimensional beam 
with the left edge fixed, see Figure~\ref{model:beam}. 
The size of the beam is 10m$\times$2m.
The domain is partitioned by triangular elements, and the equations are discretized by first order
finite elements.
The scale of this problem is $n = 46958$ with $23479$ grid nodes, where $n$ is the number of DOFs.
The first 20 smallest eigenvalues are located in [0.2526, 58.2570], and the first 100 smallest eigenvalues are located in [0.2526, 337.4505].

\subsubsection{3D hollow platform model}
\label{subsect: hollow platform model}
\begin{figure}[htbp]
\centering
\includegraphics[scale=0.22]{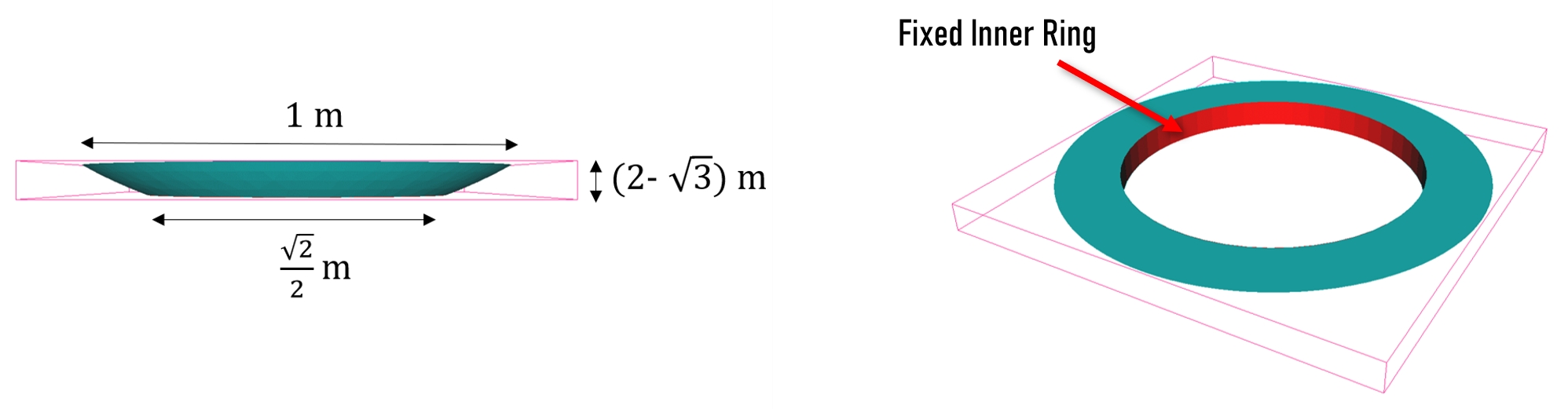}
\caption{The shape of 3D hollow platform model.}
\label{fig:hollow platform}
\end{figure}

The second model originates from cutting out a cylinder 
from a platform, see Figure~\ref{fig:hollow platform}. 
The surface of the cylinder is fixed. 
The region is partitioned into tetrahedral elements,
and second-order polynomials are used for the discretization.
The number of finite element nodes is 12387, and the dimension of the obtained GEP is 37161.
The first 20 smallest eigenvalues are located in [104.3308, 214.0376], and the first 100 smallest eigenvalues are located in [104.3308, 1506.4211].

\subsubsection{3D fish-like model}
\label{subsect: fish like model}
\begin{figure}[htbp]
\centering
\includegraphics[scale=0.24]{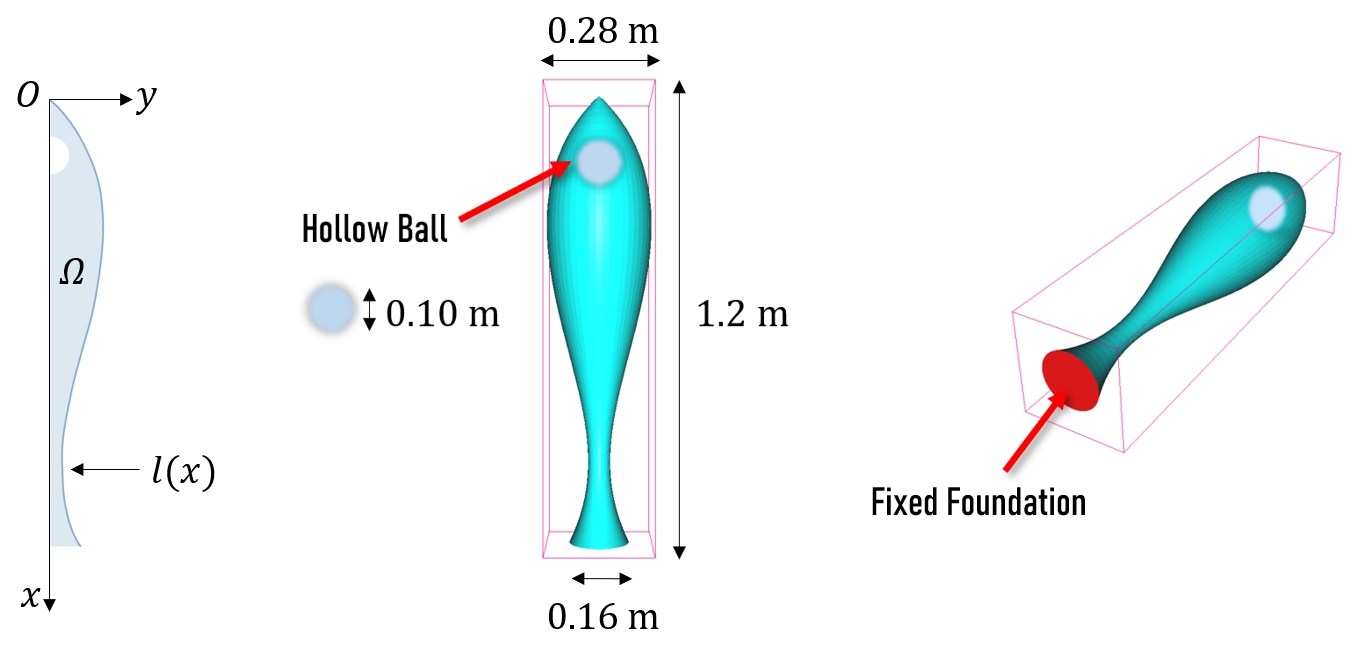}
\caption{The shape of 3D fish-like model.}
\label{fig:fish like}
\end{figure}

The third model is generated by rotating a given two dimensional region $\Omega$ along the $x$-axis, and the expression of the line $l$ is 
$l(x) = 2[ (\frac{1}{10} + (\frac{1}{3} x(x-1)^2 + \frac{1}{100} x))^{1/3}   - (\frac{1}{10})^{1/3} ]$
, see Figure~\ref{fig:fish like}. 
The obtained body looks like a fish and the lower surface is fixed.
The body is partitioned into 
tetrahedral elements, 
and a first-order finite element is used to discretize the equation.
The number of finite element node is 18416, and the dimension of the obtained GEP is 55248.
The first 20 smallest eigenvalues are located in [7.6346, 1857.5048], the first 100 smallest eigenvalues are located in [7.6346, 10238.6400].


\subsection{Numerical results}
\label{subsect:numerical-results}

\subsubsection{Setting}
\label{subsubsect: setting}
The computations were carried out in MATLAB~R2024a on a Linux compute server
with 56-core Intel Xeon Gold 6330N processors, and 1 TB of RAM.
In the following, we define $\vecr_k = \frac{(\matA - \theta_k \matB) \vecx_k}{\theta_k \Vert \matB\vecx_k \Vert_2}$
as the relative residual vector where $\theta_{k}$
being the Rayleigh quotient associated with the approximate eigenvector $\vecx_{k}$.
As stopping criteria, we use $\| \vecr_{k}\Vert < 10^{-8}$.

The details of the algorithm setting are as follows. 
Unless otherwise stated, the linear systems are solved by using the BiCGStab iterative method with an ILU preconditioner, with drop tolerance of  $10^{-4}$. 
To improve the efficiency of the ILU factorization, the matrix bandwidth is reduced by using reverse Cuthill–McKee (RCM) ordering.
The stopping criteria for the linear solver is that relative residual norm is less than $10^{-13}$ or the number of iterations has reached the maximum iteration number 1000.
In addition, we have to point out that BiCGStab might stagnate.

We get the poles and corresponding weights from the rational filter defined by $N = 4$, $\alpha = 1$ and $\beta = 0.01$.
The dimension of the search space should be larger than the number of wanted eigenvalues, so we set the number of right-hand-side terms equal to 1.2 times $\mbox{NEV}$, which is the number of wanted eigenvalues located in the interval $[a, b]$.

\subsubsection{Results}
\label{subsubsect: results}
\begin{table}[htbp]
  \centering
  \caption{Numerical results}
    \begin{threeparttable}[b]
    \begin{tabular}{c|ccccc}
    \multicolumn{1}{c}{Model} & Interval & Filter & spMV & Iter & spMV$_{\text{avg}}$ \\
    \midrule
    \multicolumn{1}{c|}{\multirow{8}[4]{*}{Beam}} & \multirow{4}[2]{*}{\shortstack{$\mbox{NEV}=20$\\(0, 58.2570]}} & Midpoint & 97530 & 13 & 7502.3\\
          &       & Gauss-Legendre & 61070 & 7 & 8724.3\\
          &       & Gauss-Chebyshev\tnote{$\dagger$} & 115120 & 12 & 9593.3 \\
          &       & SLRF  & 58042 & 11 & 5276.5\\
\cmidrule{2-6}          & \multirow{4}[2]{*}{\shortstack{$\mbox{NEV}=100$\\(0, 337.4505]}} & Midpoint & 134737 & 12 & 11228.1\\
          &       & Gauss-Legendre & 102505 & 8 & 12813.1\\
          &       & Gauss-Chebyshev\tnote{$\dagger$} & 158023 & 11 & 14365.7\\
          &       & SLRF  & 91492 & 10 & 9149.2\\
    \midrule
    \multirow{8}[4]{*}{Hollow platform} & \multirow{4}[2]{*}{\shortstack{$\mbox{NEV}=20$\\(0, 214.0376]}} & Midpoint & 17588 & 6 & 2931.3\\
          &       & Gauss-Legendre & 14891 & 5 & 2978.2\\
          &       & Gauss-Chebyshev & 18061 & 6 & 3010.2\\
          &       & SLRF  & 26105 & 9 & 2900.6\\
\cmidrule{2-6}          & \multirow{4}[2]{*}{\shortstack{$\mbox{NEV}=100$\\(0, 1506.4211]}} & Midpoint & 67944 & 8 & 8493.0\\
          &       & Gauss-Legendre & 55431 & 6 & 9238.5\\
          &       & Gauss-Chebyshev & 71339 & 7 & 10191.3\\
          &       & SLRF  & 75321 & 10 & 7532.1\\
    \midrule
    \multirow{8}[4]{*}{Fish like} & \multirow{4}[2]{*}{\shortstack{$\mbox{NEV}=20$\\(0, 1857.5048]}} & Midpoint & 60668 & 7 & 8666.9\\
          &       & Gauss-Legendre & 67595 & 6 & 11265.8\\
          &       & Gauss-Chebyshev\tnote{$\dagger$} & 95261 & 7 & 13608.7\\
          &       & SLRF  & 52283 & 9 & 5809.2\\
\cmidrule{2-6}          & \multirow{4}[2]{*}{\shortstack{$\mbox{NEV}=100$\\(0, 10238.6399]}} & Midpoint & 179471 & 14 & 12819.4\\
          &       & Gauss-Legendre & 122947 & 8 & 15368.4\\
          &       & Gauss-Chebyshev\tnote{$\dagger$} & 211025 & 12 & 17585.4\\
          &       & SLRF  & 113150 & 11 & 10286.4\\
    \bottomrule
    \end{tabular}%
    \begin{tablenotes}
        \item [$\dagger$]  \footnotesize{This symbol indicates that BiCGStab terminated due to stagnation more frequently, without reaching the required stop criterion, in this task.}
    \end{tablenotes}
    \end{threeparttable}
  \label{tab: numerical results}%
\end{table}%

The main computational cost of the rational filter methods is proportional to the number of BiCGStab iterations.
Specifically, the cost is determined by 
the sparse matrix-vector products and preconditioning operations when solving the associated multiple right-hand terms.
The reported results include the number of outer iterations \textbf{Iter},
the number of total sparse matrix-vector products \textbf{spMV} employed in the linear iterative solver, in the whole algorithm iteration.
We also record the the average number of sparse matrix-vector products \textbf{spMV$_{\text{avg}}$=spMV/Iter} for per iteration step.
In Table~\ref{tab: numerical results}, we show the numerical results when solving the first $\mbox{NEV}=20$ eigenpairs (the first sub-block) and first $\mbox{NEV}=100$ eigenpairs (the second sub-block), respectively, for the three mechanical models.

From Table~\ref{tab: numerical results}, we can see that the convergence rate of the rational filter method depends on the model problems and  the number of eigenvalues located in the specific interval. 
The Gauss-Legendre filter has the fastest convergence rate, and has the lowest computational cost for the case of the hollow platform model.
In addition to the hollow platform model, the shifted Laplace rational filter outperforms other filters in the remaining cases.
From the last column, we can see that the average cost of per iteration, \textbf{spMV$_{\text{avg}}$}, increases progressively from the newly proposed filter to the midpoint filter, Gauss Legendre filter and Gauss Chebyshev filter. This trend can be attributed to the location of the poles, which move gradually closer to the real axis, making the associated linear systems more challenging to solve.

From the numerical results above, we can also draw the conclusion that the shifted Laplace rational filter method is a good choice when solving all of the eigenpairs in the specified interval $(0, \gamma]$ for the 2D and 3D vibration models. 
\section{Conclusion}
\label{sect:conclusion}
In this paper, we develop an effective rational filter method which is suitable for computing all eigenvalues in the interval $(0,\gamma]$, for generalized eigenvalue problems arising from mechanical vibrations.

The key aspect in designing the rational filter is to restrict the poles to two straight lines, $\mathcal{L}_{\pm}: y = x(1 \pm \alpha \imath)$, with a specified slope $\alpha > 0$.
The poles and corresponding weights are then computed by solving an optimization problem based on the rational approximation of a step function.
The way to define the poles is inspired by the shifted Laplace preconditioning technique which makes the associated linear systems easier to solve by an iterative method.
Regarding the slope $\alpha$, it is important to emphasize that this parameter involves a trade-off. As $\alpha$ increases, the linear system becomes easier to solve; however, the quality of the filter deteriorates, resulting in a slower convergence rate for the eigenpairs.

Numerical results indicate that, compared to filters based on quadrature rules, this new rational filter method performs better in term of average iteration cost. 
Moreover, we observe that the difficulty of solving the associated linear systems for each pole is similar. 
This feature makes the newly proposed rational filter method highly suitable for parallel computing, as it offers excellent load-balancing capability.



 \bibliographystyle{elsarticle-num} 
 \bibliography{rationalfilter}





\end{document}